\RequirePackage{fix-cm}

\documentclass[smallextended]{svjour3}

\usepackage{etex}
\usepackage{graphicx, amssymb, amsmath, epsfig}
\usepackage{xcolor}
\usepackage{empheq}
\usepackage[T1]{fontenc}
\usepackage{pgfplots}
\pgfplotsset{compat=newest}
\usepackage{pgfplotstable} 
\usetikzlibrary{arrows,calc,backgrounds,external}
\usetikzlibrary{arrows.meta,colorbrewer}
\usepgfplotslibrary{colorbrewer}

\pgfplotsset{cycle list/Dark2-6, cycle multiindex* list={mark list*\nextlist Dark2-6\nextlist}}

\usepgfplotslibrary{external}
\tikzset{external/system call={lualatex \tikzexternalcheckshellescape --shell-escape -halt-on-error -interaction=batchmode -jobname "\image" "\texsource"}}
\tikzexternaldisable

\usepgfplotslibrary{groupplots}

\usepackage{array}
\usepackage[caption=false]{subfig}

\usepackage{datatool}
\usepackage{siunitx}
\usepackage{multirow}
\usepackage{color}
\usepackage{showlabels}
\usepackage{hyperref}
\usepackage{booktabs}
\usepackage{algorithmic,algorithm}

\usepackage{todonotes}
\usepackage{mdframed}
\usepackage{setspace}
\newcounter{mycomment}

\graphicspath{{./}{figs/}}
\newcommand{\R}{\mathbb R}

\newcommand{\cV}{\ensuremath{\mathcal V}}
\newcommand{\yl}{{\ensuremath{y^\ell}}}

\begin{document}

\title{Feedback control of parametrized PDEs via model order reduction and dynamic programming principle\thanks{The first author was supported by US Department of Energy grant number DE-SC0009324. The second and third authors would like to thank the German Research Foundation (DFG) for financial support of the project within the Cluster of Excellence in Simulation Technology (EXC 310/2) at the University of Stuttgart.}}

\titlerunning{Feedback control for parametrized PDEs via DPP}

\author{Alessandro Alla \and Bernard Haasdonk \and Andreas Schmidt}

\institute{
  A. Alla \at 
  Department of Mathematics, PUC-Rio, Rua Marquês de São Vicente, 225 Rio de Janeiro, 22453-900, Brazil \\
  \email{alla@mat.puc-rio.br}
  \and
  B. Haasdonk, A. Schmidt \at
  University of Stuttgart, Institute of Applied Analysis and Numerical Simulation, Pfaffenwaldring 57, 70569 Stuttgart, Germany \\
  \email{\{haasdonk,schmidta\}@mathematik.uni-stuttgart.de}%
}

\date{Received: date / Accepted: date}
\maketitle

\begin{abstract}
In this paper we investigate infinite horizon optimal control problems for pa\-ra\-me\-tri\-zed partial differential equations.
We are interested in feedback control via dynamic programming equations which is well-known to suffer from the curse of dimensionality.
Thus, we apply parametric model order reduction techniques to construct low-dimensional subspaces with suitable information on the control problem, where the dynamic programming equations can be approximated.
To guarantee a low number of basis functions, we combine recent basis generation methods and parameter partitioning techniques.
Furthermore, we present a novel technique to construct nonuniform grids in the reduced domain, which is based on statistical information.
Finally, we discuss numerical examples to illustrate the effectiveness of the proposed methods for PDEs in two space dimensions.
\end{abstract}

%


\section{Introduction}\label{sec:introduction}
        
Optimal control problems are challenging tasks with a huge impact in real-life applications.
The overall goal of control is to modify the behavior of dynamical systems through an external source, referred to as the control, chosen such that we are able to steer the solution trajectory to desired configurations or to achieve certain stability and performance goals.
From an application perspective, we are in particular interested in optimal controls which stabilize the system  even under perturbations.
This is a crucial point due to errors in the measurements and the inherent non-exactness of mathematical models of real-life applications.

In this work, we aim at the control of parametrized problems in feedback form, where the parameters can describe, e.g. different material parameters, geometry modifications or model uncertainties.
Usually, one is interested in solving control problems for many different parameters, e.g. in parameter studies, Monte-Carlo simulations or real-time parameter updates. This is often referred to as ``multi-query'' scenarios.

A general framework for feedback control has been introduced by Bellman in \cite{Bellman57} in the 50s via the {\em Dynamic Programming Principle} (DPP) which provides an efficient tool for the computation of the so called value function, which is an important ingredient for feedback control.
This approach is rather general and includes different optimal control problems such as, e.g. the minimum time problem, and the discounted infinite horizon control problem.
However, the method requires the solution of a nonlinear partial differential equation (PDE), e.g. the Hamilton--Jacobi--Bellman (HJB) equation, whose dimension corresponds to the dimension of the underlying control problem (see e.g. \cite{BCD97}).
Due to the nonlinearity of the HJB equation, it is usually not possible to derive analytical solutions.
Thus, it is crucial to investigate numerical algorithms to build approximations of the value functions.
Unfortunately, classical numerical methods suffer from the so-called {\em curse of dimensionality}.
Although theoretical results hold true in any dimension, the computational approximation constitutes the bottleneck of this approach.

Many numerical methods deal with the approximation of the solution to the HJB equation, such as finite volume, finite element and finite difference methods.
We refer to the monograph \cite{FF13} and the references therein for a complete presentation of suitable numerical methods.
Recently, new techniques such as radial basis functions (e.g. \cite{JS15}) and sparse grid methods (e.g. \cite{Garcke2017}) have been investigated for HJB equations.

In the current work, we will deal with semi-Lagrangian (SL) schemes which provide stable approximations of the value functions even for coarse discretizations. We make use of accelerated iterative schemes based on the fixed point iteration introduced in \cite{AFK16} where a smart coupling between a value iteration (see e.g. \cite[Appendix A]{BCD97}) scheme and policy iteration (see e.g. \cite{Bokanowski2009}) can drastically decrease the computational time to determine the numerical approximation.

Due to memory limitations, we are typically able to approximate HJB equations only up to a relatively low dimension of say $4-5$ dimensions, with a SL scheme, which is a big restriction in applications since the dimension of the HJB equation is the same as the dimension of the dynamical system.
For instance, semi-discretizations of PDEs lead to a very large number of ordinary differential equations (ODEs) which make this approach not feasible since the dimension can easily have $n \gg 10,000$ states or more, which would lead to memory requirements of the order of $\mathcal O(1/h^n)$, where $h$ is a discretization parameter, for example the grid width of a uniform grid.

One way to overcome these difficulties for PDE-related applications is to apply model order reduction (MOR) to the dynamical systems in the first place.
MOR methods (see e.g. \cite{BGW15} and the references therein) are (typically) projection based methods that have been successfully applied to different problems such as optimization and many-query problems to reduce the number of degrees of freedom of the problem and to obtain surrogate models that represent the full-dimensional and expensive model accurately. 
Although a detailed description of model reduction techniques goes beyond the scope of this work, we want to mention proper orthogonal decomposition (POD, see \cite{V13}) and balanced truncation (BT, see \cite{A05}) as two of the most popular techniques for the reduction of dynamical systems.
POD is a rather general method, which is based on a Galerkin projection of the nonlinear dynamical system onto a space whose basis functions are built upon snapshots of the system, whereas BT is based on a Petrov-Galerkin projection, where the basis functions are obtained by solving two Lyapunov equations. The reduced basis (RB) approach deals with parametric problems based on greedy algorithms (see e.g. \cite{Haasdonk13,PR07}). In this work we mainly focus on the latter approach.

The POD method has been coupled with the HJB equations in the pioneering paper studies \cite{AK01,KVX04,KX05} to compute feedback controls for high-dimensional problems for both linear and nonlinear problems.
Other features of the method have been investigated such as a-priori error estimates \cite{AFV17} and the chattering of the feedback control \cite{AFK16}.
Other model reduction methods have been coupled with the HJB approach, such as BT \cite{KK14} and more recently a comparison of reduced order modeling (ROM) techniques has been conducted in \cite{ASH2017}. 
Other approaches for the control of PDEs via the DPP deals with sparse grids for linear problems (see e.g. \cite{Garcke2017}) and spectral elements for unconstrained controls (see e.g. \cite{KK2017}).
In contrast, the SL method is rather general and includes control constraints and nonlinear dynamical systems.
As already mentioned, the focus of this work is the computation of feedback control functions for {\em parametrized PDEs} via the HJB equation coupled with MOR techniques. To the best of our knowledge this approach has not been investigated yet.

In this paper, we propose a complete workflow for the coupling of nonlinear feedback control via HJB equations and MOR. Starting from a general problem formulation, we first make use of recent ideas (\cite{SH18,ASH2017}) to project the control problem onto low-dimensional subspaces. Faced with parameter-dependent problems and with the requirement of very low-dimensional subspaces, we employ adaptive parameter partitioning techniques to reach spaces of dimension, say, maximum $5$.
For the actual numerical approximation, we employ the SL scheme for which a grid in the reduced space is required.
To this end we introduce a novel idea based on statistical assumptions on the high-dimensional system that enables data-driven approximation of the relevant part of the reduced space which is then covered by a grid.
Finally, an efficient offline/online splitting is introduced to enhance and accelerate the overall procedure.
In particular, we take advantage of the so-called value iteration (VI) scheme to precompute the value function in the barycenter of each parameter subregion in the offline phase and then switch, in the online phase, to the policy iteration (PI) method using the precomputed information on the value function as initial guess. This turns out to be a very efficient method as discussed in the numerical tests. 

To summarize, the novelties in this paper are: (i) the presence of parameters for nonlinear feedback control problems, (ii) the use of basis functions which does not depend on any control input, (iii) an automatic way to generate the domain for the reduced HJB equation, (iv) an efficient offline/online scheme, (v) numerical tests for two dimensional nonlinear equations and (vi) any initial condition.
To set the paper into perspective we recall the DPP approach and its numerical approximation in Section \ref{sec2}.
Section \ref{sec:mor_hjb} focuses on MOR for the HJB equation and all the building blocks for our approach. Finally, numerical experiments are presented in Section \ref{sec:test} with focus on the control of two-dimensional unsteady PDEs. Conclusions and future directions are discussed in Section \ref{sec:conc}.

\section{Numerical Methods for Dynamic Programming Equations}
\label{sec2}
In this section we recall the basic results for the numerical approximation of the HJB equations, additional details can be found in, e.g. \cite{BCD97} and \cite{FF13}.
Consider a continuous-time, parametric optimal control problem of the form:
\begin{gather}\label{dyn1}
  \begin{aligned}
  \min_{u\,\in\,\mathcal{U}}J_x(u;\mu),\quad\text{with}\quad J_x(u;\mu):=\int_0^{\infty} g(y(s),u(s);\mu)\, e^{-\lambda s}\,ds\\
  \qquad\text{subject to}\quad \dot y(t;\mu)=f(y(t),u(t);\mu),\quad y(0; \mu)=x,
  \end{aligned}
\end{gather}
with system dynamics $y(t;\mu)$ in $\mathbb{R}^n$ for $t\geq0$, an initial state $x\in\R^n$ and a control signal $u\in\mathcal{U}$ with
\begin{align*}
\mathcal{U}\equiv\{ u:[0,\infty)\rightarrow U, \text{measurable}\},
\end{align*}
where $U$ is a compact subset of $\mathbb{R}^m$ of admissible control values and $\lambda>0$ is the discount factor.
We consider the dynamics and the cost functional to be parametrized by a parameter vector $\mu\in\mathcal{P}\subset{\R^q}$, where $\mathcal P$ is a bounded set of admissible parameters.
The following statements and definitions are to be understood to hold for any $\mu\in\mathcal{P}$.\\
The functions $g(\cdot,\cdot;\mu)$ and $f(\cdot,\cdot;\mu)$ are assumed to be  Lipschitz-continuous functions in the first two variables.
Under rather general assumptions, the existence and uniqueness of solutions to the optimal control systems are guaranteed (see e.g. \cite{BCD97}).
A crucial tool in feedback control is the value function, which provides the minimum value of the cost functional at each point in the state space $x \in \mathbb R^n.$ For parametric problems we define it as
\begin{equation}\label{eq:VF}
v:\R^n \times \mathcal{P} \rightarrow \R,\quad v(x;\mu):=\inf_{u\in\mathcal U }J_x(u;\mu),
\end{equation}
and its characterization through the DPP for $\tau > 0$
\begin{equation}\label{dpp}
v(x;\mu)=\inf_{u\in\mathcal U}\left\{\int_0^\tau g(y_x(t,u;\mu),u(t);\mu) \mathrm e^{-\lambda t} \, \mathrm dt+v(y_x(\tau,u;\mu),u;\mu)\mathrm e^{-\lambda \tau}\right\},
\end{equation}
where we denote by $y_x(t,u;\mu)$ the dynamics of the system at time $t$ for the control signal $u\in \mathcal U$ and parameter $\mu$, starting at the initial condition $y(0;\mu)=x$. We note that we use the subscript $x$ in $y_x$ whenever we want to emphasize the dependence on the initial condition $x$. 
The above characterization can, under certain regularity assumptions on the value function, be used to derive the HJB equation: 
\begin{equation}\label{HJB}
\lambda v(x;\mu) + \sup_{u\in U}\{- f(x,u;\mu)\cdot \nabla v(x;\mu)-g(x,u;\mu)\}=0, \quad x\in\R^n,
\end{equation}
where $\nabla$ denotes the gradient with respect to $x$ from which the value function can be computed as the unique viscosity solution.
The knowledge of the value function allows the computation of the feedback control as follows:
\begin{equation*}
  u^*(x;\mu) = \min_{u\in U} \{f(x,u;\mu)\cdot \nabla  v(x; \mu) + g(x,u;\mu)\}. 
 \end{equation*}

Next, we derive a numerical scheme to approximate the value function $v(x;\mu)$.
For that purpose, we apply an SL scheme to Equation~\eqref{HJB} and thus first choose a bounded domain $\Omega\subset\R^n$ which we then discretize by a finite set of points $\Xi=\{x_i\}_{i\in J}$ with $J\coloneqq\{1,\ldots, N_G\}$ and $N_G = |\Xi|$.
We address the choice of the domain $\Omega$ and its discretization in Section \ref{sec:redom}.
Typically, $\Xi$ is a grid in $n$ dimensions and the number of grid nodes $N_G$ grows exponentially with the dimension $n$ such that already coarse discretizations lead to numbers that easily exceed the memory capacities of modern computers.
This again highlights the need for MOR techniques for high-dimensional problems. 
We construct a fully-discrete SL scheme for the approximate value function which follows from the DPP after temporal discretization of the ODEs for $y$ and a rectangular quadrature rule for the cost functional
\begin{equation}\label{HJBh} 
V(x_i;\mu)=\min_{u\,\in U}\{\mathrm e^{-\lambda \Delta t}I_1[V](x_i+\Delta t\Phi(x_i,u;\Delta t,\mu))+ \Delta t\,g(x_i,u;\mu)\},\, i\in J
\end{equation}
Here $V(x_i;\mu)$ is the approximate value for $v(x_i;\mu)$ for the nodes of the grid $\Xi$, the constant $\Delta t>0$ denotes the time-step that is used for the temporal discretization and $\Phi$ is the increment function and includes, for instance, implicit or explicit Euler schemes. 
Here $I_1[V]$ denotes a first-order interpolant of the discrete value function $V$, e.g. a piece-wise multi-linear interpolation.
This is necessary, because the point $x_i+\Delta t \Phi(x_i,u;\Delta t,\mu)$ is usually not a node of the state space grid.
Finally, let us point out how the increment function $\Phi$, introduced in Equation~\eqref{HJBh}, looks when an explicit Euler scheme is performed:
 \begin{equation}\label{incr_aff}
 \Phi(x,u;\Delta t, \mu)= f(x,u;\mu).
 \end{equation}
We refer the reader to \cite{FF13} for specific details concerning the SL schemes for HJB equations and convergence results in $L^\infty$ valid in any dimension.

We are able to approximate a solution to Equation \eqref{HJB} only up to a few dimensions, using the SL discretization with an efficient iterative solver for \eqref{HJBh}.
The simplest algorithm is based on a fixed point iteration of the value function, also called {\em value iteration} (VI):
\begin{align*}
&[V^{(j+1)}(\mu)]_i=S([V^{(j)}(\mu)]_i),\mbox{ for } j=0,1,\ldots\\
&[S(V)]_i\equiv\min_{u\,\in\,U}\{\mathrm e^{-\lambda \Delta t}I_1[V](x_i+\Delta t\Phi(x_i,u;\Delta t,\mu))+\Delta t\,g(x_i,u;\mu)\}\quad i\in J.
\end{align*}
Here we collect the nodal values in vectors $V^{(j)}(\mu) \in \mathbb R^{N_G}$, meaning $[V^{(j)}(\mu)]_i \approx V(x_i;\mu)$ where again the subscript indicates the index $i\in J$ and the superscript $j$ denotes the iteration index.
Convergence is guaranteed for any initial guess $V^{(0)}\in\mathbb R^{N_G}$ since the operator $S: \mathbb{R}^{N_G}\rightarrow\mathbb{R}^{N_G}$ is a contraction mapping (see e.g. \cite{Falcone1987}).
Although being simple and reliable, this algorithm is computationally demanding and slow when fine grids are considered.\\
A more efficient formulation is the so-called {\em policy iteration algorithm} (PI, see e.g. \cite{Howard60,SR04}), which starting from an initial guess $u^{(0)} \in U^{N_G}$ of the control at every node, performs the following iterative procedure for $i\in J$
\begin{align}\label{eq_PI}
[V^{(j)}]_i&=e^{-\lambda \Delta t}I_1[V^{(j)}](x_i+\Delta t\Phi(x_i,u^{(j)};\Delta t,\mu))+\Delta t\,g(x_i,u^{(j)};\mu)\,, \\
[u^{(j+1)}]_i&=\underset{u\,\in\,U}{\operatorname{argmin}}\{e^{-\lambda \Delta t}I_1[V^{(j)}](x_i+\Delta t \Phi(x_i,u;\Delta t,\mu))+\Delta\, t g(x_i,u;\mu)\}.
\end{align}
In the first step of \eqref{eq_PI} the PI method consists of a linear system solve since the control $u^{(j)}$ is fixed and we do not have to compute the minimization problem.
Then, the control is updated according to the value function computed in the previous step. We iterate this process until we get the desired accuracy of the value function.
It is well-known (see e.g. \cite{Bokanowski2009,SR04}) that the PI algorithm has quadratic convergence provided a good initial guess.
This point is very delicate since it requires to know a reasonable approximation of the value function.
To solve this problem we utilize an acceleration mechanism based on a VI solution on a coarse grid, which is used to generate an initial guess for PI on the fine grid, see also Section~\ref{sec:off-onl}.
This idea is based on the fact that VI generates a fast error decay when applied over coarse meshes for any initial guess.
Thus, we obtain an initial guess close to the solution and we can switch to the PI method over a fine grid, which then converges fast.
Therefore, the algorithm is a way to enhance PI with both efficiency and robustness features.
We refer to \cite{AFK15} for a detailed description of the algorithm. Finally, we note that in both algorithms we penalize the value function outside of the numerical domain to impose state constraints boundary conditions.

The main advantage of the DPP approach presented in this section is the possibility to have a synthesis of feedback controls:
Once the value function is computed, the approximated optimal control for a point $x\in\R^{n}$ in the state space is obtained by:
\begin{align*}
u^*(x)=\underset{u\,\in\,U}{\operatorname{argmin}}\{e^{-\lambda \Delta t}I_1[V](x+\Delta t \Phi(x, u;\Delta t,\mu))+\Delta\, t g(x,u;\mu)\}.
\end{align*}
For the implementation of this feedback control a direct search for the minimum can be performed if $U$ contains a finite number of control values. 
In other scenarios an efficient minimization algorithm can be employed.


\section{Framework for Parametric HJB Equations and Feedback Control}
\label{sec:mor_hjb}
\begin{figure}[htb]
  \centering
  \includegraphics[scale=.45]{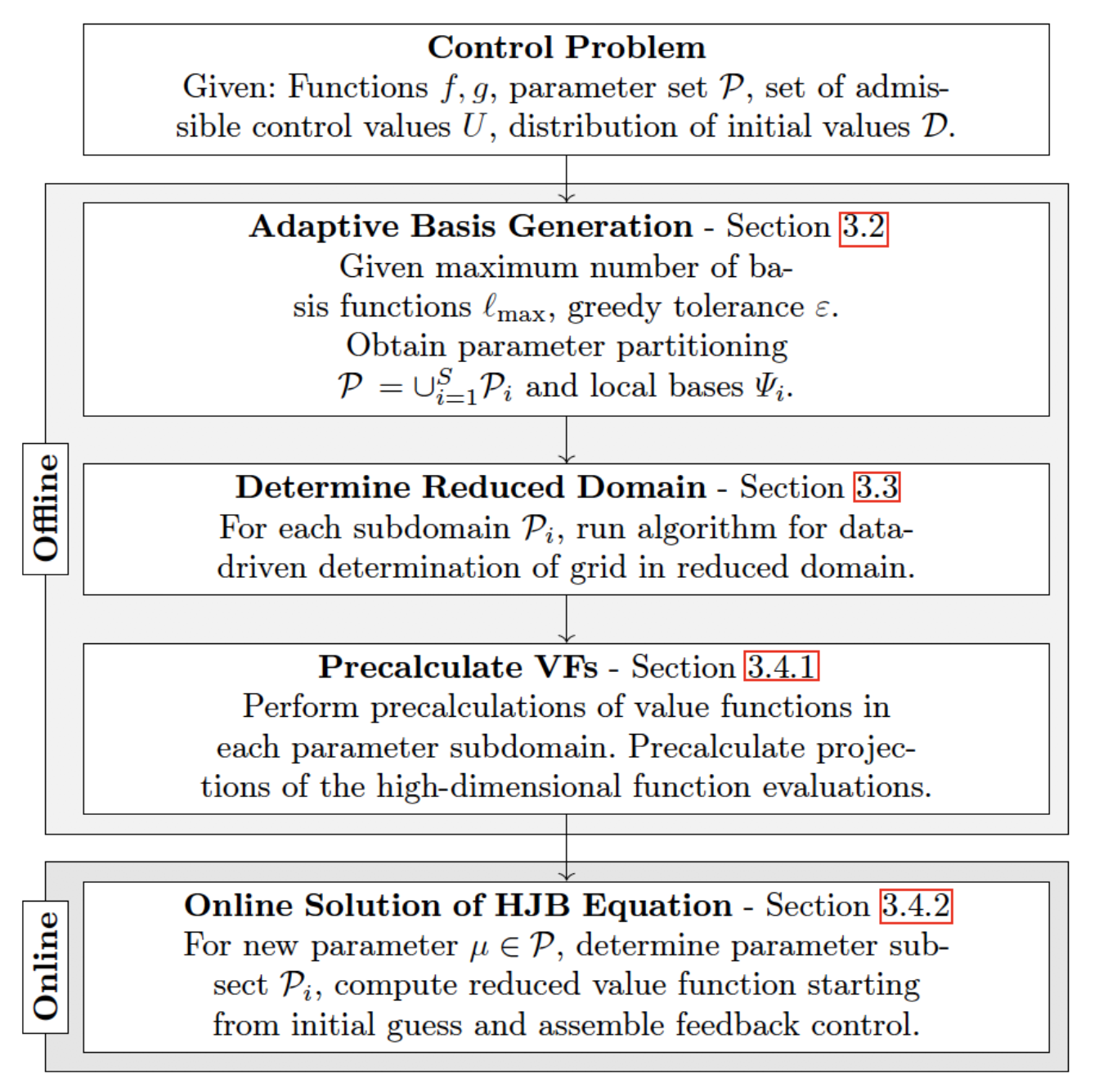}
  \caption{The workflow from a given parametric control problem towards the online approximation of the value function and feedback control.}
  \label{fig:procedure}
\end{figure}

We now provide details about the use of MOR in the context of DPP in Section~\ref{sec:projection}.
In particular, we propose a complete and automatic strategy to deal with nonlinear parametric feedback control problems via MOR and DPP.
From a computational point of view, we also show how the procedure can be implemented efficiently by employing an offline/online splitting of the whole workflow.
The overall picture of the procedure is summarized in Figure~\ref{fig:procedure}.

The first step according to Figure~\ref{fig:procedure} is an adaptive basis generation.
By this we aim at constructing suitable low-dimensional subspaces onto which we project the high-dimensional control problem.
Details about this step are explained in Section~\ref{sec:basis}.
For the numerical approximation of the value function of the so-obtained low-dimensional control problems we perform an SL scheme.
To this end we have to prescribe a finite computational domain on which the SL scheme can be applied.
To obtain knowledge about this domain we propose to make use of a data-driven method to gather statistical information about the distribution of the reduced coordinates, see Section~\ref{sec:redom}.
The actual calculation of the value function is then performed in two steps:
In a precalculation step we employ a VI scheme to get coarse approximations of the value functions in each parameter subregion, see Section~\ref{sec:offline}.
Furthermore, we precompute several expensive quantities and store them for a later reuse.
This concludes the offline step.
Online, given a parameter from the parameter domain, we refine the initial guess of the value function from the offline phase with a PI algorithm, see Section~\ref{sec:online}.
Based on the refined value function we then define the feedback control that is used to control the full dynamical system.
In what follows, we discuss in detail each of the ingredients in Figure~\ref{fig:procedure}.

\subsection{Projection-Based Approximation of the HJB Equation}\label{sec:projection}
The focus of this section is to establish a coupling between MOR and the HJB approach, as initially proposed in~\cite{KVX04} for the non-parametric case.
The need of ROM is crucial when dealing with high-dimensional problems, such as discretized PDEs, since the curse of dimensionality prohibits a direct solution of the HJB equations in higher dimensions.
We apply model reduction for the dynamical system to obtain a reduced system whose dimension is then feasible to approximate the HJB equation. 
The ROM is based on projecting the nonlinear dynamics onto a low $\ell$-dimensional subspace $\cV\subset \R^n$ that contains the relevant information about the dynamics $y(t;\mu)$.
We equip the space $\cV$ with an orthonormal basis, given by the columns in the matrix $\Psi \in \R^{n\times \ell}$, which will be specified in the following section. We then approximate the full state vector by the linear combination of basis vectors, i.e. $y(t;\mu) \approx \Psi \yl(t;\mu)$ where $\yl:[0,\infty)\rightarrow \R^\ell$ are the so called reduced coordinates.
Plugging this ansatz into the dynamical system and requiring a Galerkin condition, we obtain an ODE-system of dimension $\ell$
\begin{align}
  \begin{aligned}\label{eq:reduced_sys}
  \dot{y}^\ell(t;\mu)&= \Psi^Tf(\Psi \yl(t),u(t);\mu),\quad t>0,\cr
  \yl(0;\mu)&=\Psi^T x.
  \end{aligned}
\end{align}
The procedure presented above is a generic framework for MOR for dynamical systems. We note that it is possible to extend the whole procedure performing a Petrov-Galerkin projection and we refer to \cite{ASH2017} for a detailed description of the method.
It is clear, that the quality of the approximation highly depends on the reduced space $\cV$, i.e. on the chosen basis $\Psi$.
In particular, the dynamical system \eqref{eq:reduced_sys} should capture enough information to allow for accurate approximations of the closed-loop behaviour for any desired initial state. 

To ease the notations we introduce abbreviations for the reduced quantities.
The initial value will be denoted by $x^\ell := \Psi^T x$ whereas the projected dynamical system and the reduced running cost function are given as
\begin{align*}
    f^\ell(y^\ell(t;\mu),u(t);\mu)&:= \Psi^T f(\Psi y^\ell(t;\mu), u(t); \mu), \\
    g^\ell(y^\ell(t;\mu), u(t);\mu)&:=g(\Psi y^\ell(t;\mu),u(t);\mu).
\end{align*}
In the general projection framework above, we define the optimal control problem for the projected system:
\begin{gather}
  \inf_{u \in \mathcal U} J^\ell_{x^\ell} (u;\mu) := \inf_{u \in \mathcal U} \int_0^\infty g^\ell(\yl(t;\mu),u(t);\mu)\mathrm e^{-\lambda t}\; d t,\\
  \begin{aligned}
    \text{s.t.}\quad  \dot{y}^\ell(t;\mu) &= f^\ell(y^\ell (t;\mu), u(t);\mu),  \quad t>0, \\
    \yl(0;\mu) &= x^\ell.
  \end{aligned}\label{eq:redocp}
\end{gather}
As in the full-dimensional case, we define the value function for the reduced system
$$v^\ell(x^\ell;\mu): = \inf_{u \in \mathcal U} J^\ell_{x^\ell}(u;\mu),$$
which satisfies the reduced HJB equation which is now $\ell$-dimensional and feasible for numerical treatment as long as the dimension is sufficiently small, e.g. $\ell\leq 5$
\begin{equation}\label{redHJB}
\lambda v^\ell(x^\ell;\mu) + \sup_{u\in U}\{- f^\ell(x^{\ell},u;\mu)\cdot \nabla_{x^\ell} v^\ell(x^\ell;\mu)- g^\ell(x^\ell,u;\mu)\}=0,\quad\forall x^\ell\in\R^\ell.
\end{equation}
We note that the reduced HJB equation \eqref{redHJB} will admit a unique viscosity solution as the original problem since it is obtained by orthogonal projection. We refer to \cite{AFV17} for more details.

The overall idea is now to replace the high-dimensional value function $v(x;\mu)$ by its reduced counterpart $v^\ell(x^\ell;\mu)$.
Furthermore, we make use of the reduced value function and define the following approximated feedback law, which is essentially the control law from the full-dimensional system, where the value function is replaced by the low-dimensional approximation:
\begin{equation}\label{eq:feed}
   u^\ell(x) = \min_{u\in U} \{f(x,u;\mu)\cdot \nabla  v^\ell(x^\ell;\mu) + g(x,u;\mu)\}. 
 \end{equation}
In Section~\ref{sec:test}, we will show the quality of the feedback control, when applied to the full-dimensional system.

To obtain computationally efficient schemes, we assume that the dynamics has a linear dependence with respect to the control $u$
\begin{equation}\label{dyn:dec}
    f(y,u;\mu) = f^y(y;\mu) + f^u(y;\mu) u,
\end{equation}
and that the cost functional has a quadratic form
\begin{equation}\label{quad_cost}
g(y,u;\mu) = y^T Q(\mu) y + u^T R(\mu) u,
\end{equation}
where $Q(\mu)\in\mathbb R^{n\times n}$ is symmetric positive semidefinite and $R(\mu)\in\mathbb R^{p\times p}$ is symmetric positive definite.

\subsection{Basis Generation Algorithm}\label{sec:basis}
Let us now provide more information about the computation of the basis functions $\Psi$.
Since our numerical schemes are li\-mi\-ted to a very low number of basis functions, the quality of the basis is of utmost importance.
In~\cite{ASH2017} a comparison for different basis generation techniques in the context of feedback control via the HJB equation is performed in the non-parametric context. It turns out that classical but straightforwad approaches such as POD or BT do not necessarily yield satisfying results since the focus of those methods is on providing surrogate models for the dynamics and not for feedback control. A different approach that is based on the explicit form of the value function in the linear case is introduced in \cite{ASH2017}.

Finding the basis functions is of course rather hard for arbitrary nonlinear control problems. Unlike for linear problems with quadratic cost functionals where the value function can be computed explicitly by solving an algebraic Riccati equation (ARE), the value function for nonlinear problems is in general not known analytically.
However, we can always obtain local information about the basis by linearizing the control problem around a constant point of interest $(\bar y, \bar u)$:
\begin{align*}
  f(y,u; \mu) \approx f_y(\bar y, \bar u;\mu)(y - \bar y) + f_u(\bar y, \bar u; \mu)(u-\bar u).
\end{align*}
In the sequel, we will typically choose $(\bar y, \bar u) = 0$, since we are interested in steering the system to the origin and hence can write the linearized state equation as
\begin{align*}
  \dot y = A(\mu)y + B(\mu) u,
\end{align*}
with matrices $A(\mu) \coloneqq f_y(\bar y, \bar u; \mu) \in \R^{n\times n}$, $B(\mu) \coloneqq f_u(\bar y, \bar u;\mu) \in \mathbb R^{n\times p}$. 

Solving the linearized optimal control problem can be easily done by means of the associated discounted ARE for the positive semi-definite and stabilizing solution $P(\mu)\in \R^{n\times n}$ 
\begin{align}\label{eq:ricc}
 \mathcal R(P(\mu)) := \left(A(\mu)-\dfrac{\lambda}{2} I_n\right)^T P(\mu) + P(\mu)\left(A(\mu)-\dfrac{\lambda}{2} I_n\right) -\nonumber\\
  P(\mu)B(\mu)R(\mu)^{-1}B(\mu)^T P(\mu) + Q(\mu) = 0,
\end{align}
where $\mathcal R(P(\mu))$ is the residual of the ARE, $Q(\mu)$ and $R(\mu)$ given accordingly to \eqref{quad_cost}. We refer to \cite{FF13} for a detailed derivation of the ARE under the presence of the discount factor $\lambda$.
The idea is now to build basis functions upon the information on the value function from the linearized control problem for varying parameters, although it is just an approximation to the true and unknown value function.
For that purpose, we adopt the low-rank factor greedy (LRFG) procedure from~\cite{SH18} to our setting.
For the sake of completeness, we summarize the method in Algorithm~\ref{alg:lrfg}.
It runs in a typical greedy structure: An error indicator is minimized over a suitably large but finite training set $\mathcal P_\text{train} \subset \mathcal P$ of parameters by adding information about the worst-approximated true solution in each iteration.
In line 4 of the algorithm, only the part which is not yet captured in the basis is considered by orthogonalization and in step 5, the remaining information is compressed via an additional POD, where we prescribe a desired level $1-\varepsilon_\text{POD}$ of {\em POD-energy} for some $\varepsilon_\text{POD} \in [0,1]$ that should be captured by the basis, see e.g.~\cite{SH18,V13} for details. The algorithm is a variant of the POD-Greedy procedure, which is known to be quasi-optimal for MOR of parametric unsteady PDEs (see e.g. \cite{Haasdonk13}).
We choose the error indicator as the normalized residual norm $\Delta(\mu) := \| \mathcal R(\hat P(\mu) ) \|_F/\|Q(\mu)\|_F$ where $\hat P(\mu) \in \R^{n\times n}$ is the approximate solution to the ARE for the current basis. We note that the error indicator proposed here was certified in \cite{SH18}.
\begin{algorithm}[tbh]
  \begin{algorithmic}[1]
    \REQUIRE Parameter training set $\mathcal P_\text{train} \subset \mathcal P$, desired greedy tolerance $\varepsilon$, POD tolerance $\varepsilon_\text{POD}$, initial basis $\Psi$
    \WHILE{$\max_{\mu \in \mathcal P_\text{train}} \Delta(\mu,\Psi) > \varepsilon$}
    \STATE $\mu^\ast := \operatorname{arg\,max}_{\mu \in \mathcal P_\text{train}} \Delta(\mu)$
    \STATE Solve ARE for $P(\mu^\ast)$
    \STATE $P_\perp := (I-\Psi \Psi^T) P(\mu^\ast)$
    \STATE $\Psi := [\Psi, \operatorname{POD}(P_\perp, \varepsilon_\text{POD})]$
    \ENDWHILE
  \end{algorithmic}
  \caption{LRFG algorithm for the calculation of the projection basis.}
  \label{alg:lrfg}
\end{algorithm}
This procedure might be expensive since it requires the solution of an ARE and a subsequent SVD in each iteration.
However, by employing low-rank techniques for the solution of the large-scale AREs, both of these tasks can be sped up substantially.
We also refer the interested reader to the recent work \cite{Simoncini2016b} for an in-depth discussion of projection-based model reduction for the ARE and the link to the LQR problem.
As already mentioned, the strength of this model reduction approach relies on the fact that the basis functions contain directly information of the value function for the infinite horizon problem.
Furthermore, the whole described technique does not depend on a particular choice of the control, unlike POD. 

By applying a grid-based sche\-me for the approximation of the value function, we are restricted to a relatively low number of dimensions $\ell$ for which the procedure can be performed.
Furthermore, the presence of parameters can change the control problem significantly when going from one configuration to another.
Therefore, a basis which is able to capture information about the whole parameter domain might easily exceed the maximum possible dimension.
To overcome this problem, we apply an adaptive method introducing a partitioning of the parameter domain $\mathcal{P}$.
By running the adaptive algorithm, the parameter domain $\mathcal P$ is split into $S$ partitions $\mathcal P_i \subset\mathcal P$ such that $\mathcal{P}=\cup_{i=1}^{S} \mathcal{P}_i$ together with local bases $\Psi_i, \in \R^{n \times \ell_i}$ where $\Psi_i$ are the $\ell_i$ basis functions computed for the partition $\mathcal{P}_i$ for $i=1, \ldots, S$.

The idea behind the partitioning is that we want to deal with a prescribed maximum number $\ell_\text{max}$ of basis functions for each subregion of the parameter domain to guarantee the computational feasibility of the reduced control problem\eqref{eq:redocp} and, simultaneously guarantee a certain accuracy $\varepsilon$. 
The algorithm works as follows:
Given a partitioning $\{\mathcal P_i\}_{i=1}^S$ which initially is set to $\mathcal P$, we run the basis generation on each parameter subset $\mathcal P_i$ independently.
Two cases can occur:
Either the desired accuracy is reached within the prescribed number of basis functions $\ell_{\text{max}}$, or the error indicator/the number of basis functions is too large. 
In the latter case the parameter region $\mathcal P_i$ is refined for example by bisection and the procedure is repeated on all newly identified subregions.
The method stops when the desired accuracy is reached and the number of basis elements $\ell_\text{max}$ is not exceeded in each subdomain.
As stopping criterion we also include a maximum number of refinements since it is not always possible to reach the amount of basis functions required.
In these cases we accept a reduced basis of lower accuracy that satisfies a strict size constraint $\ell_i\leq \ell_\text{max}$. 
We refer to e.g. \cite{HDO11,EPR10a} and the references therein for more details.

\subsection{Data-Driven Approximation of the Reduced Domain}\label{sec:redom}
In this section, we provide details on the procedure to determine the domain for the approximation of the reduced HJB equation \eqref{redHJB}.
Although the reduced HJB equation is defined on the full space $\mathbb R^\ell$, for numerical reasons we have to restrict ourselves to a bounded domain $\Omega^\ell\subset\R^\ell$ and the question arises how a reasonable choice can be made.
Note that the design of the reduced domain $\Omega^\ell$  is also of great importance for the application in feedback control:
By applying the reduced-order feedback control from \eqref{eq:feed}, the projection $\Psi^T y_x(t;\mu) \in \R^\ell$ of the current state of the controlled system onto the reduced space is required and fed into the reduced value function.
Hence, in order to get accurate feedback controls, those projected vectors should be contained in the domain $\Omega^\ell$ where the reduced HJB equation is approximated. Finally, the use of a domain requires to impose boundary conditions to the reduced HJB equation. A common choice is the penalization of trajectories which exit the domain, see e.g. \cite{AFV17}.

A common approach is to choose the domain a-priori of the form $\Omega^\ell = [a_1,b_1]\times \ldots \times [a_\ell, b_\ell]$ where $a_i < b_i,\, i=1,\ldots,\ell$ are prescribed bounds.
However, it is not clear how the values $a_i$ and $b_i$ can be chosen.
In particular, in a parametric scenario, the influence of the parameter can greatly alter the dynamics and thus the projections.

In the current work we propose a novel strategy that makes use of statistical information about the full and reduced coordinates.
For that purpose, we assume that the initial values of the high-dimensional problem follow a prescribed multivariate distribution, which we abbreviate by $x \sim \mathcal D$ where $\mathcal D$ defines the probability density chosen.
This choice is motivated by the following heuristic observation: In cases where the states of the discretized system represent nodal values, e.g. in a FE scheme, the values of neighboring nodes are often of very similar magnitude.
This results from phenomena like diffusion or other types of transport of information. In other scenarios, often statistical a-priori knowledge of the states that can occur in the application are available, e.g. typical temperatures in a heat transfer application.
Note that we have to consider some assumptions on the full states $y\in\mathbb R^n$, since otherwise their projections $y^\ell = \Psi^T y$ can lie anywhere in $\mathbb R^\ell$ and the approximation may be arbitrarily bad.

\begin{algorithm}
\caption{Data-driven approximation of the numerical domains.\label{Alg:dom}}
\begin{algorithmic}[1]
  \REQUIRE Parameter partitioning $\mathcal P = \cup_{i=1}^S \mathcal P_i$ with local bases $\Psi_i$, distribution $\mathcal{D}$, time instances $T = \{t_0,\dots,t_K\}$, desired number of grid nodes $H_{i,j}$.
  \FOR{$i=1,\dots,S$}
    \STATE Choose $\mu^\ast$ from $\mathcal P_i$.
    \STATE $Y \gets [y_{\xi}(t_k;\mu^\ast)]_{\xi \in X, t_k \in T}$ (collect snapshots).
    \STATE $\tilde Y \gets \Psi_i Y$.
    \FOR{$j=1,\dots,\ell_i$}
    \STATE $h_j \gets $ approximation of distribution of $j$-th component.
    \STATE $\Phi_j \gets$ univariate grid $\{ s_1, \dots, s_{H_{i,j}} \}$ with $0\in \Phi_j$ and $\int_{s_q}^{s_{q+1}} h_j(s)\mathrm ds$ equal for $q=1,\dots,H_{i,j}-1$.
    \ENDFOR
    \STATE Build non-uniform grid $\Xi_i \coloneqq \Phi_1 \times \dots \times \Phi_{\ell_i}$.
  \ENDFOR
\end{algorithmic}
\end{algorithm}

\begin{figure}
  \centering
\includegraphics[scale=.45]{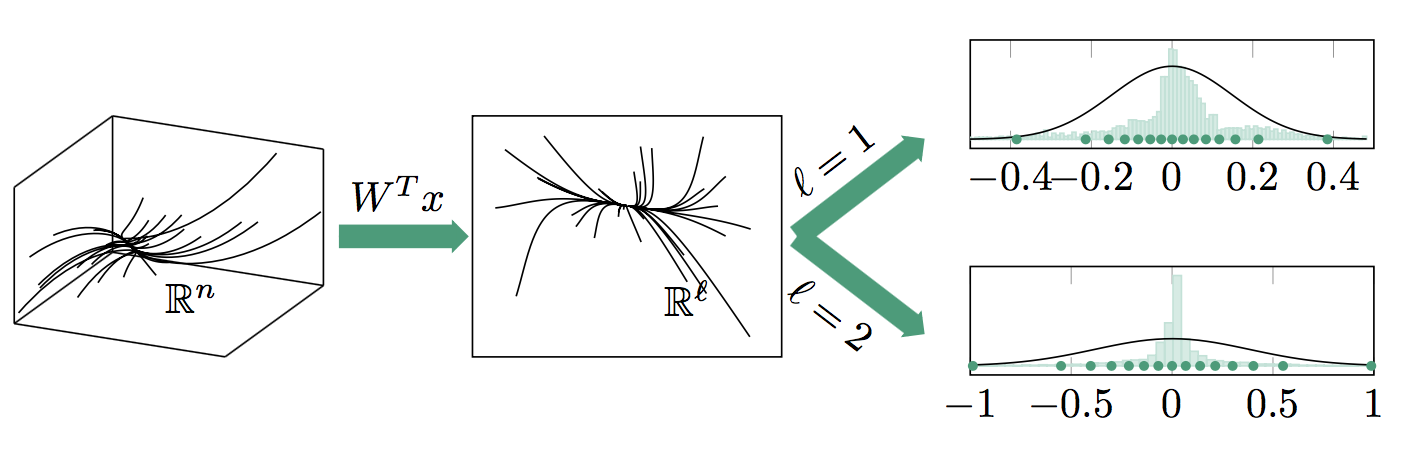}
  \caption{The procedure for obtaining the statistical distribution and estimating the univariate distributions ($N=15$ partitions in this case).}
  \label{fig:redom_procedure}
\end{figure}

The idea of the proposed algorithm is to sample solutions to the high-dimensional system for certain suitable controls and parameters and to estimate the componentwise distribution of the projected reduced vectors.
The algorithm is summarized in Algorithm~\ref{Alg:dom} and illustrated in Figure~\ref{fig:redom_procedure}.
Based on the given parameter partitioning, we loop over all $S$ parameter regions and perform the following procedure, where $i$ always denotes the index of the current parameter partition:
First, we pick a sample parameter $\mu^\ast$ from the $i$-th parameter domain, e.g. the barycenter and a set of $N_\text{train}$ initial conditions $X=\{\xi_1,\ldots, \xi_{N_\text{train}}\}$ with $\xi_k\sim\mathcal{D},$ for $k=1,\ldots, N_\text{train}$.
Then we simulate the (controlled) high-dimensional system with the parameter $\mu^\ast$ and all $\xi \in X$ and collect the solution at time instances $T$ in a snapshot matrix $Y$.
The control $u^\ast$ can for example be chosen from the linearized system for $\mu^\ast$ or simply be set to zero in case of stable systems.
We then project the snapshots onto the $\ell_i$-dimensional subspace that is spanned by the basis $\Psi_i$ and analyze the result componentwise:
To this end, we normalize the data and fit a Gaussian to the distribution of the values in the reduced coordinate.
From this we get a continuous function $h_j : \R \rightarrow [0,1]$ with $\int_\R h_j(s) \mathrm ds=1$ which we use to construct a set of grid nodes $\Phi_j$ for the component.
Given a desired odd number of grid nodes $H_{i,j}$, we enforce the area under the curve of $h_j$ to be equal between all grid nodes.
By this we ensure a distribution of the grid nodes that fits to the estimated statistical information.
The final grid $\Xi_i$ for the parameter region $\mathcal P_i$ is then defined as the cartesian product of all one-dimensional grids and consists of $|\Xi_i|=\prod_{j=1}^{N_i} H_{i,j}$ points.
A schematic drawing of the procedure for the reduction of $n=3$ to $\ell=2$ is given in Figure~\ref{fig:procedure}.

\subsection{Offline-Online Efficient Implementation of the Numerical Scheme}\label{sec:off-onl}
In this section we provide some remarks about technical features of the method to improve its computational efficiency. In particular, we will explain how to deal with an offline-online decomposition that is often used in MOR to distinguish the first phase characterized by potentially expensive computations to build a surrogate model (offline stage), which enables rapid and inexpensive simulations (online stage). 
In the current work, the role of model reduction is slightly different, in fact we are interested in reducing the dimension of the dynamical systems to decrease the complexity of the corresponding HJB equation and, therefore, to make the problem feasible.
We do not aim at real-time computations since the method will still rely on the approximation of a high-dimensional PDE.
However, by assuming a special structure of the ODE function $f$ and the running cost $g$, we can realize a speed-up since the expensive evaluations of the nonlinear function $f$ can be shifted to the offline stage. We recall that an offline/online decomposition, in this context, is new and we will show the computational benefit in Section \ref{sec:test}.
We further assume that the dynamics given in \eqref{dyn:dec} satisfies the following parameter separability assumption:
\begin{align}\label{aff_dec}
  \begin{aligned}
    f(y,u;\mu) &= f^y(y;\mu) + f^u(y;\mu) u = \sum_{q=1}^{Q_y} \Theta_q^y(\mu) f_q^y(y) + \sum_{q=1}^{Q_u} \Theta_q^u(\mu)f_q^u(y) u,
  \end{aligned}
\end{align}
where the functions  $\Theta^y_i$, $\Theta^u_j: \mathcal{P}\rightarrow\R$ for $=1,\ldots, Q_y$ and $j=1,\ldots, Q_u$ are coefficient functions depending only on the parameter $\mu$. 
This structure allows for the precomputation of most function evaluations that are needed during the online phase. 

\subsubsection{Offline Stage}\label{sec:offline}
The offline stage constitutes the building block of our approach, where most of the quantities are precomputed and stored for any parameter configuration. It basically consists of three parts:
\begin{enumerate}
\item basis generation, including the parameter partitioning, 
\item sampling of the set for the estimation of the grid in the reduced space,
\item preparatory tasks for a fast online PI, which is explained in the following.
\end{enumerate}

Given a parameter partitioning $\mathcal P_i$ for $i=1,\ldots,S$, together with corresponding grids generated by the procedure explained in Section~\ref{sec:redom}.
Let us denote the corresponding grid nodes corresponding to the $i$-th subdomain $\mathcal P_i$ as $\Xi_i$ with $|\Xi_i|=H_i$.
In order to speed up the online calculations for the PI, we make use of the special structure defined in Equation~\eqref{aff_dec} and precompute all function evaluations and their projections
\begin{align*}
  f_{q,i}^{y,\ell} := \left[ \Psi_i^T f_q^y(\Psi_i x_1), \dots, \Psi_i^T f_{q,i}^y(\Psi_i x_{H_i}) \right], \quad q = 1,\dots,Q_y, \\
  f_{q,i}^{u,\ell} := \left[ \Psi_i^T f_q^u(\Psi_i x_1), \dots, \Psi_i^T f_{q,i}^u(\Psi_i x_{H_i}) \right], \quad q = 1,\dots,Q_u,
\end{align*}
where $x_j \in \Xi_i$ for $j=1,\ldots, H_i$ and $i=1,\ldots, S$. 
Note that the pre-calculated quantities are of low-dimension and can all be precomputed once in the offline phase.
Given a parameter $\mu$, and using assumption \eqref{aff_dec} the function can be rapidly evaluated on the grid nodes for this parameter by summing up the parameter-independent quantities $f_{q,i}^{y,\ell}$ and $f_{q,i}^{u,\ell}$, weighted by the corresponding coefficient functions $\Theta_q^y(\mu)$ and $\Theta_q^u(\mu)$.

This allows us to reduce drastically the number of evaluations of the dynamical system with powerful speed up in both the VI and PI method since the computation of the value functions involves several evaluations of the dynamical systems. Similar assumptions may be posed on the cost functional. 
We finally note that \eqref{aff_dec} is absolute crucial to perform fast evaluations and can be always obtained via the EIM algorithm (see e.g. \cite{BMNP04}).


The next step concerns the precalculation of the value function via a VI scheme.
As described in Section \ref{sec2} we solve the reduced HJB equation \eqref{redHJB} performing a VI algorithm and then switching to a PI method to obtain fast convergence of the method.
We propose to use the VI method offline for some particular choices of the parameter.
In fact, since we act with a partition of the parameter domain we assume that in each subregion the dynamics will not differ significantly: since the basis generation yielded a low-dimensional basis we compute an approximation of the value function for the barycenter of each subdomain $\mathcal{P}_i$.
Therefore, we obtain accurate initial guesses for the value function and guarantee fast convergence when switching to the PI algorithm on a finer grid.
Finally, we note that at this stage we compute the finer grid, on which we later use the PI method and evaluate all the quantities independent from the parameter $\mu$ as e.g. $f_q^{y,\ell}, f_q^{u,\ell}$. 

%

\subsubsection{Online Stage}\label{sec:online}
The precomputation in the offline stage allows us to focus on the following steps in the online phase:
Given a new parameter $\mu \in \mathcal P$ we have to
\begin{enumerate}
\item identify the parameter partition $\mathcal P_i$ such that $\mu \in \mathcal P_i$,
\item calculate an accurate approximation for the reduced value function $v^\ell(\cdot;\mu)$,
\item define the feedback control $u^\ell(x)$ according to equation~\eqref{eq:feed}.
\end{enumerate}

The first step is trivial and the last step can be readily performed once the approximation $v^\ell(\cdot;\mu)$ is available.
For the second step we run a PI algorithm starting from the initial guess for the value function, which was calculated during the offline phase.
Note that at this point we can make use of the precalculated function evaluations on the grid $\Xi_i$ to speed up the calculation significantly.
By doing this, the overall complexity does not depend on the high-dimension $n$ but only on the reduced dimension $\ell_i$ and the number of grid points in $\Xi_i$.
%
We compute the reduced value function $v^\ell(\cdot;\mu)$ satisfying \eqref{redHJB} at each grid point $x^{\ell_i} \in\Xi_i$ .
\begin{equation*}
u^{\ell_i}(x;\mu):=\arg\min_{u\in U}\big\{f(x,u;\mu)\cdot \nabla_{x^{\ell_i}} v^{\ell_i}(x^\ell;\mu)+g(x,u;\mu)\big\}.
\end{equation*}

We note that here we replace the high-dimensional value functional with the reduced approximation whereas the dynamics $f$ and the cost functional $g$ are kept high-dimensional due to the fact that the basis functions better describe the value function rather than the dynamics. This strategy turns out to be more stable than using the reduced functions $f^\ell$ and $g^\ell$. 

\section{Numerical Tests}\label{sec:test}

We now present three examples of optimal feedback control problems, which demonstrate the efficiency of our proposed method.
The first example models a control problem for a linear advection-diffusion equation where the true optimal feedback control and the true value function can be computed by means of the ARE.
Thus, we are able to compare the numerical approximation obtained from our approach to the true solution.
The second example is a two-dimensional semi-linear heat equation with a cubic nonlinearity which presents an unstable equilibrium around the origin.
The control objective will be the stabilization around this point.
The last example deals with a coupled viscous Burgers system, introducing many layers of additional complexity since the dimension of the control space is two as well as the number of outputs.
Furthermore, the equations for this scenario are described by two coupled nonlinear PDEs.
The aim of the second and third example is to show that nonlinear feedback control is more efficient than a LQR controller based on the linearization of the problem which is then plugged into the nonlinear model under consideration.

\begin{figure}
\includegraphics[scale=.45]{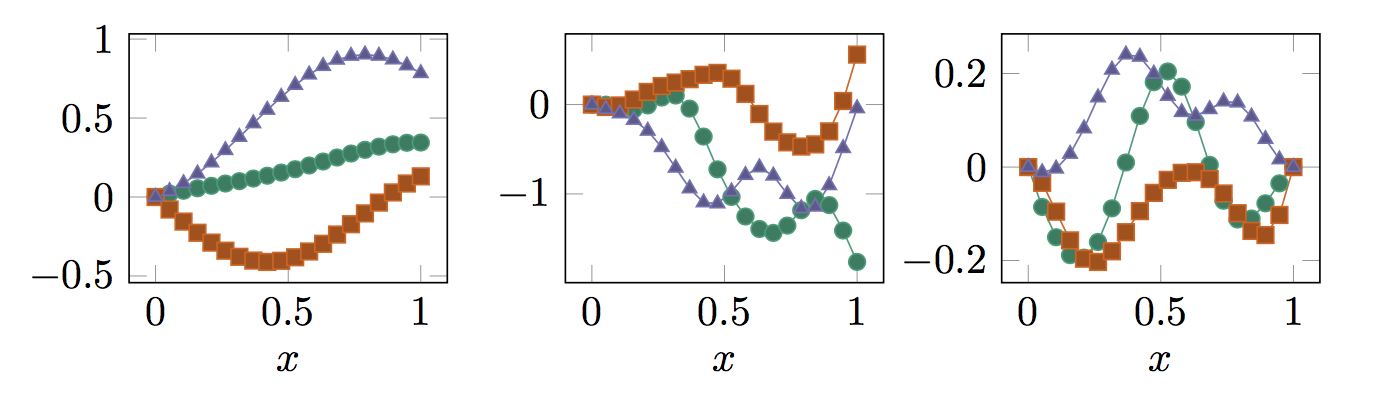}
  \caption{Three random initial conditions chosen by the distribution $\mathcal N$ for different $\gamma$ and $b(\cdot)$ choices: $\gamma=1,b(x) = x$ (left); $\gamma=20,b(x) = x$ (middle); $\gamma=10, b(x) = -\frac{1}{4} + (x-\frac{1}{2})^2$.\label{fig:example_covariance}}
\end{figure}

To apply the workflow we have to assume certain statistical properties of the high-dimensional solution.
Since our examples stem from semi-discretized PDEs, we can define those properties based on the nodal values of the discretization.
To this end let $N_1,\dots,N_n \in \R^d$ be the coordinates of the nodes in either the FE mesh or in the FD discretization.
We then define the Gaussian distribution $\mathcal D \coloneqq \mathcal N(\nu, \Sigma)$ with mean $\nu \in \R^n$ and positive-definite covariance matrix $\Sigma \in \R^{n \times n}$.
To model the relationship between the nodes we define the entries in the matrix $\Sigma$ as
\begin{align}
  \Sigma_{i,j} = c \cdot b(N_i)b(N_j) \mathrm e^{-\gamma \|N_i-N_j\|^2}, \quad i,j = 1,\dots,n, \label{eq:gaussian_cov}
\end{align}
for $\gamma, c >0$.
By changing $\gamma$ we adjust the weight of nodes that are close/far to each other.
Furthermore, through the function $b:\R^d \rightarrow [0,\infty)$ we get the possibility to put different weights on nodes that are, for example, close to the boundary.
By this we can incorporate zero boundary conditions.
Figure~\ref{fig:example_covariance} shows three random vectors drawn according to the distribution $\mathcal D$ for a discretization of the interval $[0,1$] into $n=20$ nodes with different values for $\gamma$, different weight functions $b$ and $\nu=0$.
The example shows the great flexibility in the modeling of the high-dimensional states.

\subsection{Test 1: Two-Dimensional Linear Advection-Diffusion}\label{sec:num_lin}
\begin{figure}
  \centering 
\includegraphics[scale=.45]{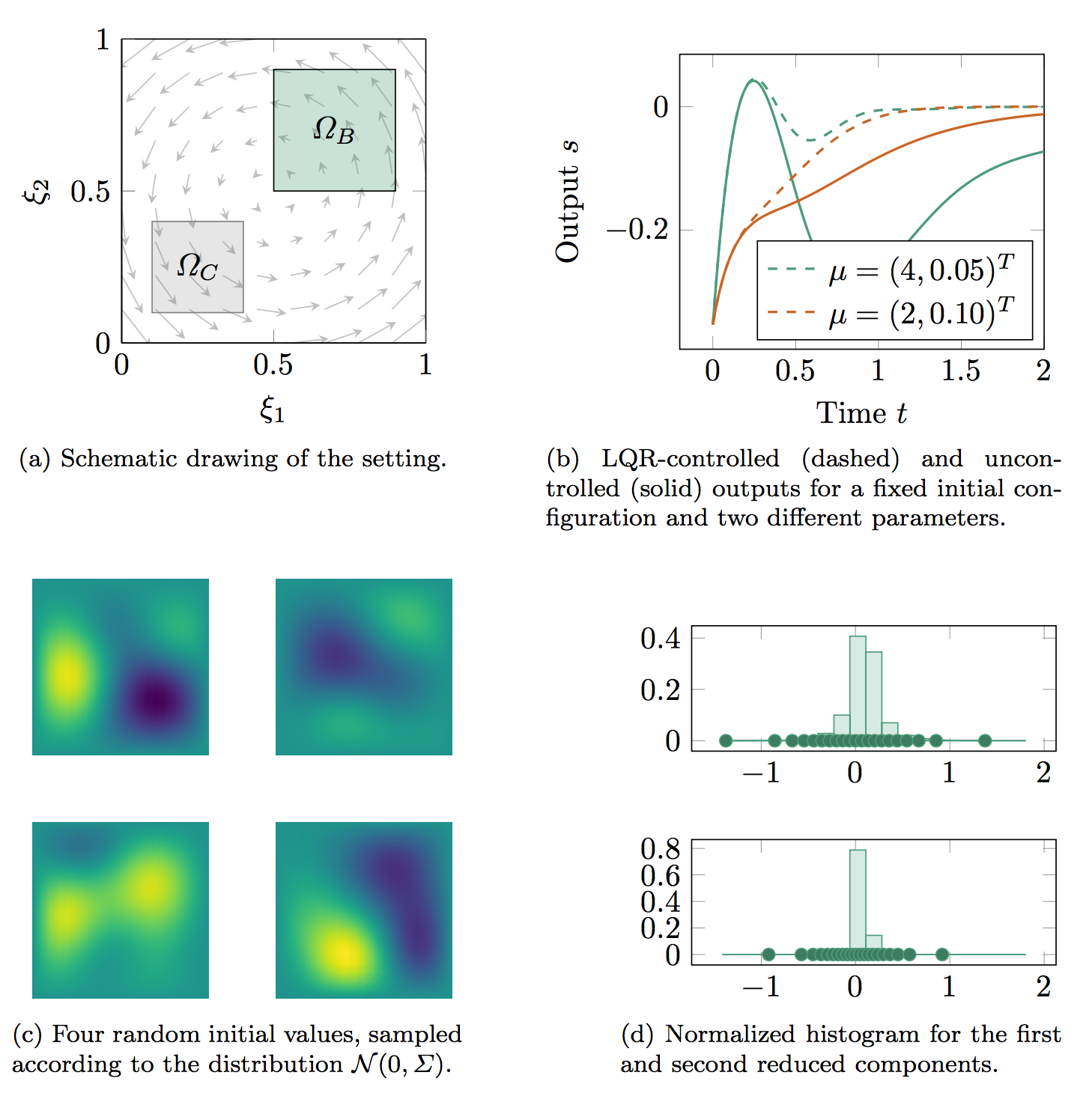}
  \caption{Linear setting.}
  \label{fig4}
\end{figure}

The first test problem considers an optimal control scenario for a two-dimensional linear advection-diffusion equation on the domain $\Omega := (0,1)^2$.
The parametric PDE for this example is given by
\begin{align*}
 \partial_t w(t,\xi;\mu) - \mu_\text{diff}\Delta w(t,\xi;\mu) + a(\xi;\mu) \cdot \nabla w(t,\xi;\mu) &= \mathrm 1_{\Omega_B}(\xi) u(t), \quad \xi \in \Omega, t\geq 0,\\
 w(0,\xi;\mu)&=w_0(\xi;\mu).
\end{align*}
The velocity is defined as the divergence-free field at $\xi=(\xi_1, \xi_2)\in\Omega$ as
$$a(\xi;\mu) = \mu_{\text{adv}}\cdot(-(\xi_2 - 0.5), (\xi_1-0.5))^T,$$
which induces a counterclockwise flow in the solution with velocity $\mu_{\text{adv}}$ as shown in Figure \ref{fig4}.
The indictor function $ \mathrm 1_{\Omega_B}(\xi)$ maps the scalar control $u(t)$ onto the domain $\Omega_B := [0.5,0.9]^2$.
The parameters in this example are $(\mu_\text{diff},\mu_\text{adv})^T \in \mathcal P:=[0.05,0.1]\times[2,4]$.
In order to set up the control problem we define the standard quadratic cost functional 
\begin{align*}
  J_{w_0}(u;\mu) := \int_0^\infty \mathrm( 10 \,s(t;\mu)^2 + 10^{-2} u(t)^2 )\mathrm e^{-\lambda t}\mathrm dt
  \intertext{with}
  s(t;\mu) := \frac{1}{|\Omega_{C}|}\int_{\Omega_{C}} w(t,\xi;\mu) \,\mathrm d\xi,\, t\geq 0,
\end{align*}
where $\lambda=10^{-3}$ and $\Omega_C := [0.1,0.4]^2$. Here, we introduce the quantity of interest $s(t;\mu)$ that depends on the solution of the PDE $w(\cdot,\cdot;\mu)$. 
Figure \ref{fig4} shows the outputs for two different configurations for the controlled and uncontrolled problem. We see how different parameters lead to different outputs. 
We discretize the control problem in space by using linear finite elements on a uniform triangular grid, resulting in an $n=676$ dimensional LTI system with the scalar discretized output $z$
\begin{align}\label{test1:disc}
  E\dot y = (\mu_\text{diff}A_\text{diff} - \mu_\text{adv}A_\text{adv} - A_\text{dirichlet}) y + B u,\quad z = C y,\quad y(0) = x.
\end{align}
We note that Equation \eqref{test1:disc} fits into the abstract setting shown in \eqref{dyn1} with $Q=10 C^T C, R=10^{-2}$ in \eqref{quad_cost} and that $A_{\text{diff}}, A_\text{adv},$ and $A_\text{dirichlet}$ are the discretization of $\Delta w$, $\nabla w$ and the boundary conditions, respectively. 
The temporal discretization is carried out with an implicit Euler scheme with step size $\Delta t=10^{-2}$.
In what follows, we compare the controlled dynamics with both the HJB and LQR approach. 
For this setting it is known that the true value function is given as $v(y;\mu) = y(0)^TP(\mu)y(0)$ and the true optimal control takes the form $u(t) = -K(\mu) y(t)$ with the feedback gain matrix $K(\mu):=R^{-1}B^TP(\mu)$ that depends on the solution $P(\mu)\in \mathbb R^{n\times n}$ of the ARE \eqref{eq:ricc}.
To apply the HJB approach we restrict the control values to the finite set of points $U := \{ u^3\,|\, u=-2+i\Delta u , i=0,\dots,109 \}$ where $\Delta u = \frac{4}{109}$, which provides enough information to capture the LQR control values sufficiently accurate.

For all simulations, we choose initial values that are sampled via a Gaussian distribution with covariance matrix \eqref{eq:gaussian_cov} with zero mean and the choices with $c=10^{-3}, \gamma=2$ and $b(y) := (-4(y_1 - 0.5)^2+1)(-4(y_2 - 0.5)^2+1)$ for the node coordinates $y:=(y_1,y_2)^T\in\R^2$ of the FE mesh.
In Figure~\ref{fig4} we provide a plot of four random initial vectors, drawn by using the distribution $\mathcal N(0,\Sigma)$. 
According to our proposed method, we first run the adaptive LRFG algorithm to produce a partition, $\mathcal P = \cup_{i=1}^S \mathcal P_i$ of the parameter space and corresponding local bases. 
The basis generation was performed with a desired tolerance of $\varepsilon=0.9$, maximum basis size $5$ and maximum refinement level $3$, resulting in a grid as indicated in Figure~\ref{fig4}.
The next step consists of gathering statistical information for building the grids in each parameter subregion $\mathcal P_i$. For that purpose we run $100$ uncontrolled simulations and collect the reduced state vectors.
Figure \ref{fig4} shows the distribution of the uncontrolled solutions for the parameter subdomain for $\mu=(\mu_\text{diff},\mu_\text{adv})^T$. We can see that the distribution allows to compute a non-uniform grid which is finer where the distribution is higher.
The domain for the reduced HJB equation is computed as discussed in Section \ref{sec:redom}.

Let us first investigate the performance of the HJB approach compared to the true LQR for a fixed parameter. For this purpose, we choose the test parameter $\mu^\ast = (3,0.08)^T$ which leads to a non-trivial configuration due to the large advection and small diffusion.
We calculate a fixed basis by solving the ARE for this parameter and using the first $\ell=\{1,2,3,4,5\}$ left singular vectors as basis elements, see also Section~\ref{sec:basis}.
We show the results of the controlled problem and compare the error in the costs of the full-dimensional system steered with the LQR control and the approximated HJB control in Table~\ref{tbl:lin_fixed_par}. 
For that purpose we pick initial vectors from a test set 
$$\mathcal X := \{x^i \in \mathbb R^n | x^i \sim \mathcal N(0,\Sigma), 1\leq i \leq N_\text{test}:=100\}$$
and run the full dimensional simulations with the controls obtained by a reduced LQR controller, e.g. the optimal control law of the reduced order LTI system, and the approximated HJB controls. 
As error measurement, we define the mean relative error in the costs of the approximated controlled systems $J_x(u;\mu)$ compared to the true LQR cost $J^\text{LQR}_x(u;\mu)$ for the initial state $x$ as
$$\operatorname{mean}\limits_{x \in \mathcal X} \frac{|J_x(u, \mu) - J^\text{LQR}_x(u,\mu)|}{|J^\text{LQR}_x(u,\mu)|}.$$
The first column in Table~\ref{tbl:lin_fixed_par} represents the dimension of the reduced problem. The second column is the error when using the LQR controller which is obtained by solving the ARE for the reduced order system. The third, fourth, fifth and sixth column show the error between the true value function and value function computed by the HJB-approach with 31 points in each dimension (third and fourth column) and only 11 points (fifth and sixth column). As one can see, the reduced problem of dimension 1 and 2 is not stable for both, the LQR and HJB approach. Then increasing the dimension $\ell$, the error decays as expected.
It is also possible to see that our proposed approach for the discretization of the reduced domain performs better than the equidistant grid. This is a consequence of a finer grid around the point of interest. We also note that our results are very close to the LQR which we consider optimal here. However, it is hard to make a fair comparison between the LQR and the HJB approach because their settings are not the same e.g. the control space and the numerical domain.
Furthermore, the table shows the quality of the basis functions for a strong advection dominated problem.
\begin{table}
     \centering
  \begin{tabular}{cccccc}
       \toprule
       &  & \multicolumn{2}{c}{HJB with $31$ points} & \multicolumn{2}{c}{HJB with $11$ points}\\
      $  \ell $ & LQR & Equi. & Non-Equi. & Equi. & Non-Equi.\\
      \midrule
   $1$ & $1.18 \cdot 10^0$ & $1.98 \cdot 10^3$ & $1.10\cdot 10^3$ & $2.27\cdot 10^3$ & $1.72\cdot 10^3$\\
   $2$ & $4.34 \cdot 10^1$ & $4.69 \cdot 10^2$ & $3.92\cdot 10^2$ & $3.86\cdot 10^3$ & $5.87\cdot 10^2$\\
   $3$ & $2.53 \cdot 10^{-1}$ & $1.79 \cdot 10^{-1}$ & $1.42\cdot 10^{-1}$ & $2.36\cdot 10^{-2}$ & $1.78\cdot 10^{-1}$\\
   $4$ & $2.12 \cdot 10^{-2}$ & $5.25 \cdot 10^{-2}$ & $2.87\cdot 10^{-2}$ & $1.02\cdot 10^{-1}$ & $6.90\cdot 10^{-2}$\\
      $5$ & $2.76 \cdot 10^{-3}$ & $4.23 \cdot 10^{-2}$ & $2.07\cdot 10^{-2}$ & $9.27\cdot 10^{-1}$ & $5.10\cdot 10^{-2}$\\
\bottomrule
\end{tabular}
   \caption{Test 1: Approximation error for the fixed parameter $\mu^\ast=(3,0.08)^T$. }
   \label{tbl:lin_fixed_par}
\end{table}

\begin{figure}
  \centering
\includegraphics[scale=.45]{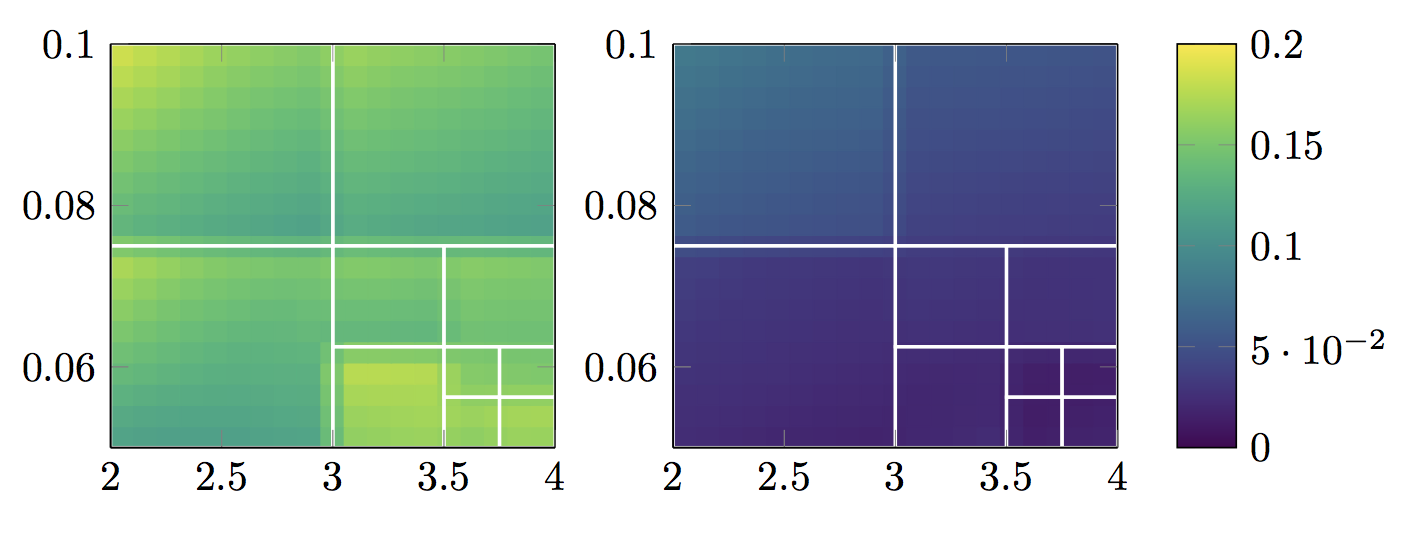}
  \caption{Test 1: Relative error of the value function computed with VI algorithm in the barycenter of each parameter subregion in the offline stage (left), relative error of the value function computed with a PI algorithm over the whole parameter space in the online stage (right).}
  \label{fig:linear_parameter_space_evaluation}
\end{figure}

Then, let us draw our attention to the parametrized problem. As discussed in Section \ref{sec:basis} we have to use an adaptive strategy not to exceed a certain number of basis functions to be able to solve the reduced HJB equation. Figure \ref{fig:linear_parameter_space_evaluation} shows how the algorithm identifies subregions in the parameter space. It is somehow intuitive that advection dominated problems need more information on the basis functions and therefore further refinements towards higher $\mu_\text{adv}$ and lower $\mu_\text{diff}$.
Here, we decide to use $\ell=4$ basis functions in each subregion.
Finally, we show the error over the whole parameter space in Figure~\ref{fig:linear_parameter_space_evaluation}.
In the left panel one can see the error of the value function computed with an offline VI algorithm only in the barycenter of each subregion. In the right panel, we see how the PI algorithm improves the approximation of the value functions. We also note that the VI is computed on a very coarse grid with $9$ points in each dimension, whereas the PI algorithm is computed with $25$ points. This plot also shows the benefit of the offline/online decomposition in terms of accuracy of the value function.

\subsection{Test 2: Nonlinear Unstable 2D Heat Equation}
\begin{figure}
  \centering
\includegraphics[scale=.45]{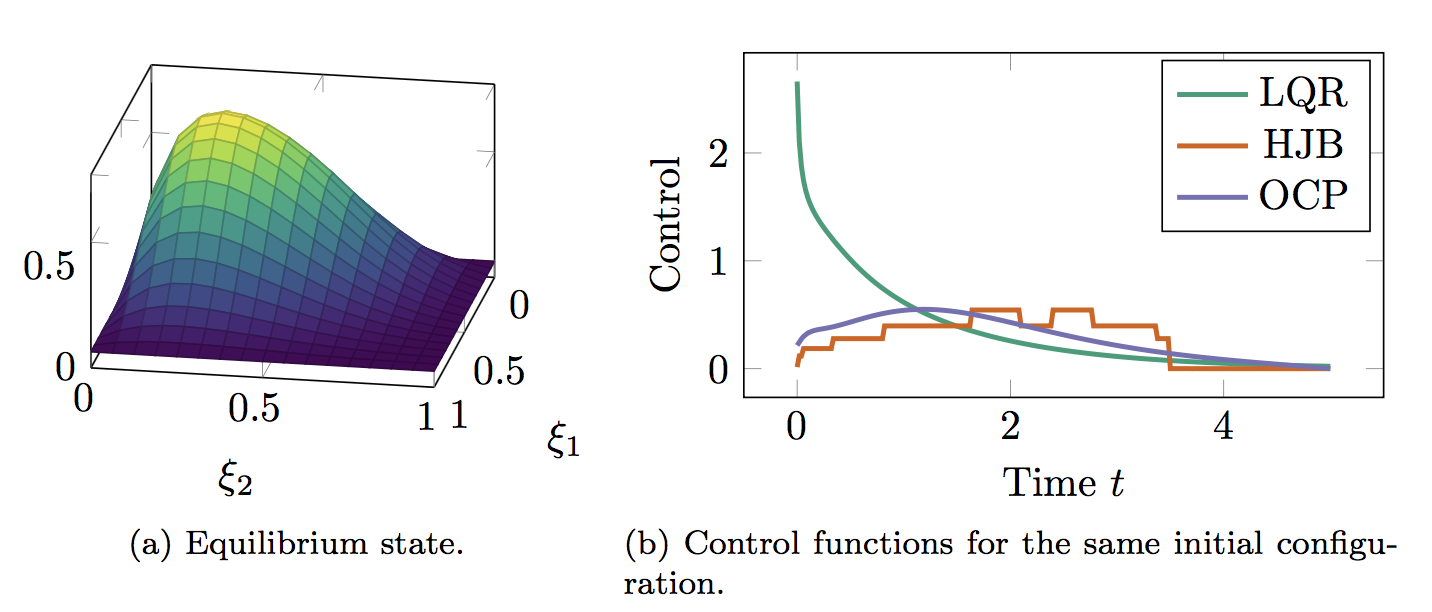}
  \caption{Test 2: Equilibrium state and control values for $\mu=7$.}
  \label{fig2}
\end{figure}

The second test problem deals with the control of a two-dimensional semi-linear advection-diffusion equation with a cubic nonlinearity on the domain $\Omega := (0,1)^2$.
The parametric PDE is defined as
\begin{align}\label{eq:nonlinheat_pde}
\begin{aligned}
 \partial_t w(t,\xi;\mu) - \mathcal L w(t,\xi;\mu) + \mu \left(w(t,\xi;\mu)-w(t,\xi;\mu)^3\right) = \mathrm 1_{\Omega_B}(\xi) u(t),\\
  \xi \in \Omega, t\geq 0,
 \end{aligned}
\end{align}
where the linear operator $\mathcal L w := 0.2\Delta w - \nabla\cdot w$ describes the diffusion and advection part.
We impose homogeneous Dirichlet boundary conditions on all boundaries and define the distributed control input via the indicator function $\mathrm 1_{\Omega_B}(\xi)$ on the domain $\Omega_B := [0.2,0.6]^2$. 
The parameter $\mu$ directly influences the strength of the cubic nonlinearity and takes values in the set $ \mathcal P:=[2,7]$. We consider $Q=10, R=1$ and the discount factor is $\lambda = 10^{-3}$ in \eqref{quad_cost}.
The problem is spatially discretized by using a finite difference scheme on a uniform grid with $n=361$ nodes.
The resulting system takes the form $\dot y = A y + \mu F(y) + B u$, which fits to our assumption about the offline/online splitting for the online PI.
Here, the nonlinearity $F(y)$ is the component-wise evaluation of the cubic nonlinearity, e.g. $(F(y))_i = (y_i - y_i^3)$.
The temporal discretization is performed by applying an explicit Euler scheme with step size $\Delta t=10^{-3}$.
Equation~\eqref{eq:nonlinheat_pde} has 3 equilibria, where $w=0$ is unstable.
The uncontrolled dynamics reach either the stable equilibrium depicted in Figure~\ref{fig2} or the one which has the same structure but opposite sign.
The control goal is to steer the solution to the unstable origin and keep it there.
We also note that this particular example is not stable under finite-time open-loop control, since it is impossible to reach exactly zero and any small deviation will lead to instabilities.
We also mention that a model predictive control approach can be applied as an alternative to our proposed method, see e.g.~\cite{Gruene2011}. 
For the sake of completeness we also show the control input in Figure \ref{fig2} computed with LQR, HJB and an open-loop approach. We note that the HJB approach agrees with the open-loop solution. However, the latter method is computationally very expensive since we have to chose a long horizon to approximate the infinite horizon problem and this might lead to unstable solutions. For these reasons, in what follows, we only compare our results with the LQR routine. We recall that the LQR controller is computed from the linearized system and plugged into the nonlinear model.
In this example, we again make use of a Gaussian distribution with covariance matrix \eqref{eq:gaussian_cov} with $b$ chosen as in the linear and with $\gamma=5$ and $c = 0.45$.

We apply the procedure proposed in this paper and therefore start with the basis generation. 
The linearization around the origin of the state equation yields an LTI system of the form $\dot y= (A + \mu I_n) y + Bu$.
Applying the adaptive LRFG algorithm for this example is not trivial because the solutions to the ARE do not have a low-rank structure but full numerical rank.
A heuristic explanation for this is that we measure the full state in the cost functional instead of an output of interest. 
However, we successfully apply the algorithm and prescribe a desired level of refinement and run the algorithm.
The algorithm then refines uniformly over the parameter space up to the prescribed level.
We run this procedure with maximum refinement levels $1$ and $3$, resulting in two and $16$ partitions.

\begin{figure}
  \centering
\includegraphics[scale=.45]{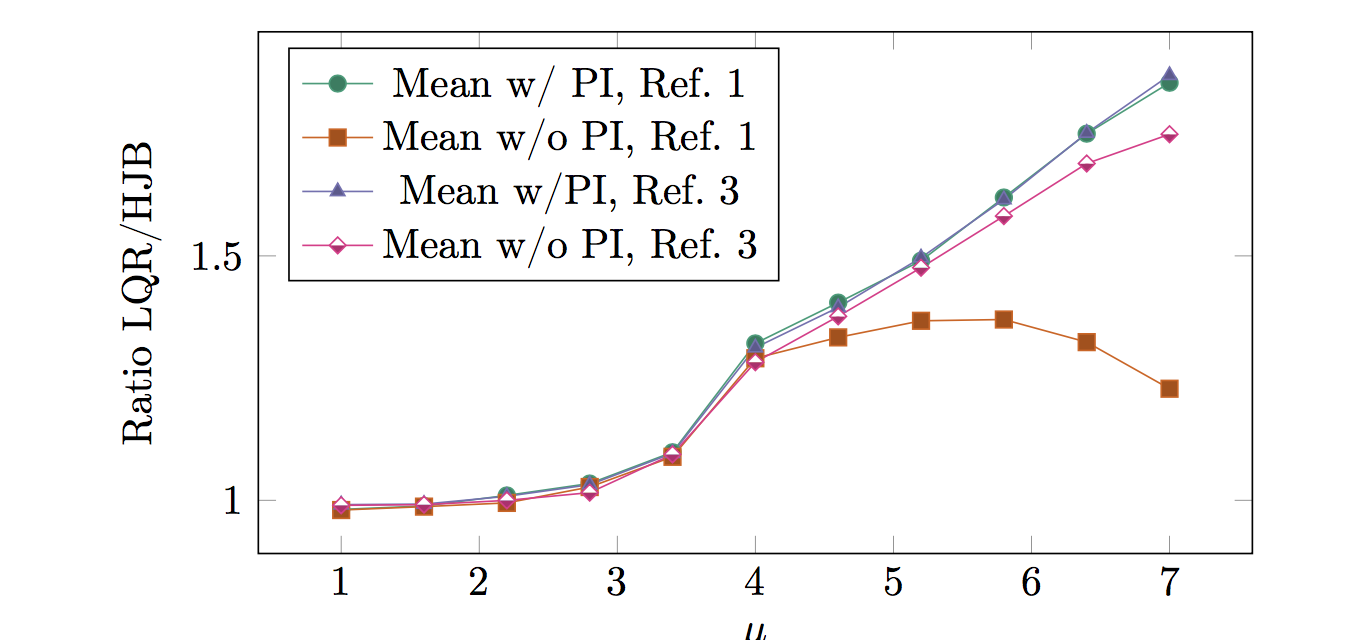}
  \caption{Test 2: Plot showing the improvement of the HJB approach vs. a classical LQR approach.  \label{fig:nonlinheat_improvement}}
\end{figure}

In Figure~\ref{fig:nonlinheat_improvement} we show the average ratio $\operatorname{mean}_{x \in \mathcal X}J^{LQR}_x / J^{HJB}_x$ for $100$ random samples to demostrate the improvement of the HJB approach over using a classical LQR controller obtained from the linearized system.
As one can see, for small values of the parameter $\mu$, our results are very close to the LQR setting due to a small contribution of the nonlinear term.
However, when increasing $\mu$, we can observe a huge improvement with the HJB approach.
Furthermore, we see how the refined parameter partitioning influences the accuracy of our approximations.
Even though both refinement levels yield the same improvement over the LQR-controlled simulation after running an online PI, we see that we are able to reach the same quality of approximation by just using the value functions calculated offline in the barycenters for the third refinement case.
From a computational point of view, we note that the online PI starting from the coarse approximations with one refinement level requires much more iterations to converge.

\subsection{Test 3: Burgers' System}
\begin{figure}[htbp]
  \centering
 \includegraphics[scale=.5]{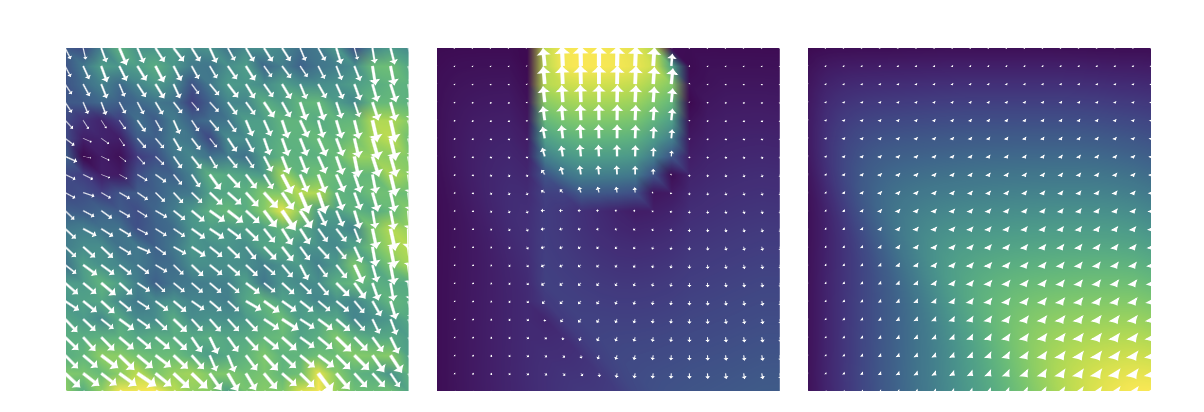}
   \caption{Test 3: Example of the dynamical behaviour of the Burgers' equation. From left to right: Initial state, controlled state and uncontrolled state at $t=2.5$.}
  \label{fig3:burg}
\end{figure}

\begin{figure}[htbp] 
\includegraphics[scale=.45]{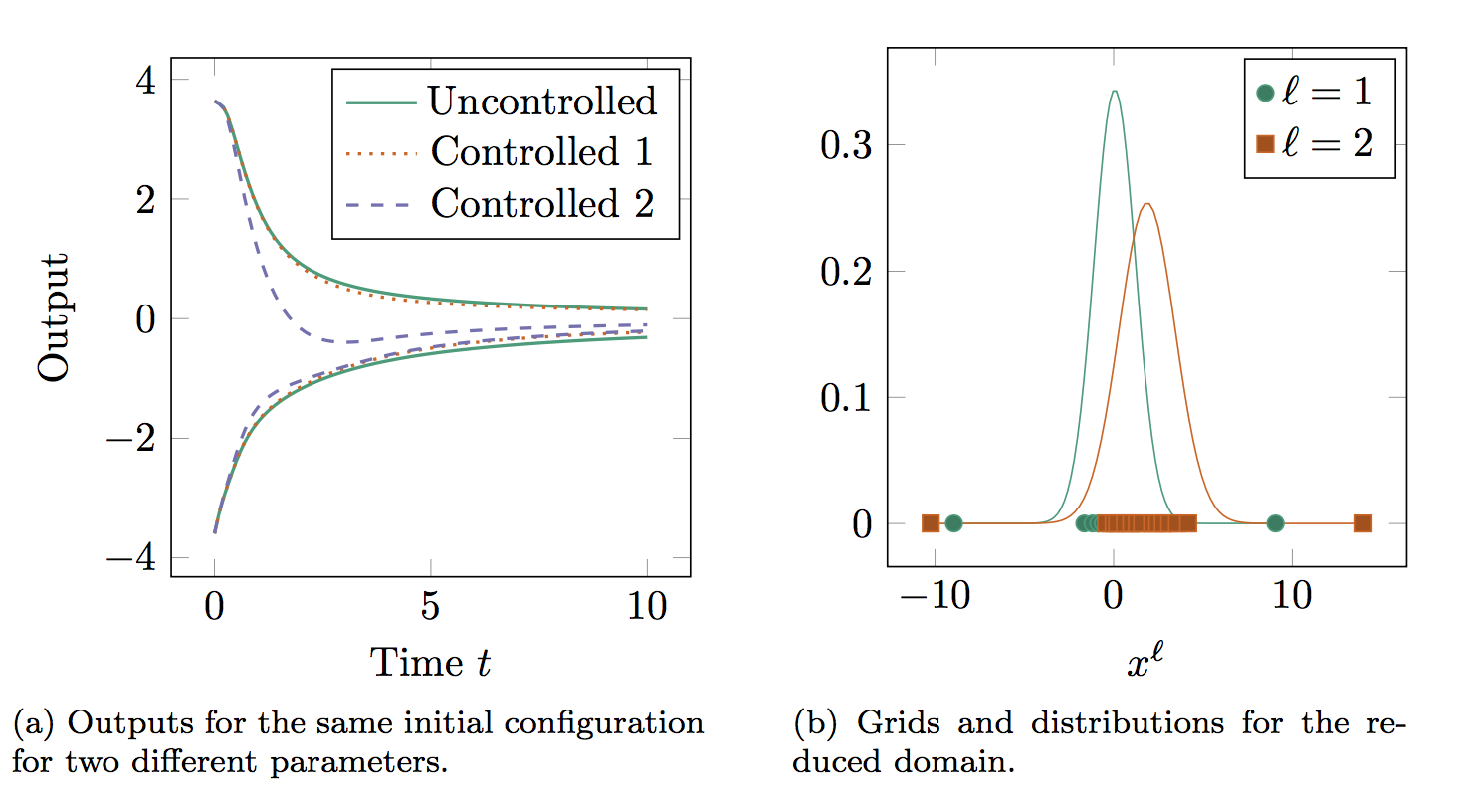}
  \caption{Test 3: Outputs of the controlled and uncontrolled example and the distributions and grids that were generated.}
  \label{fig3}
\end{figure}

We consider the following two-dimensional coupled Burgers equations for $\xi \in \Omega \coloneqq (0,1)^2$ and $t \geq 0$
\begin{gather}\label{eq:burg}
  \begin{aligned}
  \partial_t w(t,\xi;\mu) - \sigma \Delta w(t,\xi;\mu) + (w(t,\xi;\mu) \cdot \nabla )w(t,\xi;\mu) &= \mathrm 1_{\Omega_B}(\xi)u(t),\quad \\
  w(0,\xi;\mu) &= w_0(\xi;\mu),
  \end{aligned}
\end{gather}
for the unknown function $w(t,\xi;\mu) = (w_1(t,\xi;\mu), w_2(t,\xi;\mu))^T \in \mathbb R^2$.
We impose homogeneous Dirichlet boundary conditions on all boundaries and choose a low diffusion constant $\sigma = 10^{-4}$.
The indicator function $1_{\Omega_B(\xi)}$ maps the two control functions $u(t) = (u_1(t),u_2(t))^T$ onto the subdomain which is given, component-wise, by the ball of radius $0.2$ centered in $(0.5,0.25)^T$, i.e. $\Omega_B = B_{0.2}( (0.5, 0.25)^T )$.
We consider two partial measurements $s_1(t)$, $s_2(t)$ from the system:
\begin{align}\label{eq:burgers-output}
  s(t; \mu) := \left(\mu_1 \int_{\Omega_C} w_1(t,\xi;\mu) \mathrm d\xi, \quad \mu_2 \int_{\Omega_C} w_2(t,\xi) \mathrm d\xi\right), \quad t\geq 0,
\end{align}
which are the average velocities of the flow in $\xi_1$ and $\xi_2$ direction, measured on the subdomain $\Omega_C := B_{0.2}( (0.5,0.25)^T )$ and the parameters $\mu_1,\mu_2 \in \mathcal{P}:=[0.01, 5]$ determine weights on the individual flow components.
The cost functional for the PDE control problem is given by:
$$ \int_0^\infty \mathrm e^{-\lambda t} (\|s(t;\mu)\|^2 + \|u(t)\|^2) \mathrm dt, $$
with the discount factor $\lambda = 10^{-4}$.
We discretize the system \eqref{eq:burg} by a finite difference scheme with an upwind flux for the convection term which leads to a system of ODEs of dimension $n = 800$ 
 \begin{align*}
  \dot y(t) = f(y(t),u(t)) = \sigma Ay(t) + Bu(t) + F(x(t)),\quad z(t;\mu) = C(\mu) y(t),
\end{align*}
where the discretized output $z(t)$ stems from a discretization of Equation~\eqref{eq:burgers-output} by a rectangular quadrature rule.
Since we consider two inputs and two outputs, the dimension of the matrices are given as $B\in\R^{n\times 2}$ and $C(\mu) \in \R^{2\times n}$.
As control space we choose $U = \bar U^2$ where $\bar U = \{ u^3 | u = -3 + 0.1875i, i = 0,\dots,32 \} \}$.
The temporal discretization is carried out with an explicit Euler scheme with time step $\Delta t=5\cdot10^{-3}$.
For this example we again have to specify a distribution of initial values of interest.
In both components, we pick Gaussian distributions with a covariance matrix given by equation \eqref{eq:gaussian_cov} with $c = 0.2$, $b=1$ and $\gamma = 1$.
In the second component we furthermore set the mean of the distribution to $-1$ yielding flows which are mostly directed from top to bottom. 
Figure~\ref{fig3:burg} shows an example simulation, both the uncontrolled and HJB-controlled case.
In Figure \ref{fig3} we show two outputs for two different parameter configuration.

The ROM is built upon the LRFG algorithm on the linearized equation from \eqref{eq:burg}. The algorithm yields a surrogate model of dimension $\ell=2$ and does not perform a parameter partitioning.
We note that the linearization of~\eqref{eq:burg} around the origin leads to the heat equation which is possible to reduce with a few basis functions. For the numerical domain of the reduced HJB equation, we compute 100 uncontrolled test simulations where the initial values were chosen according to the distribution $\mathcal{D}$.
Then, in each dimension the reduced domain is discretized with 15 points according to the distributions that were estimated from the reduced coordinates of the test simulations, see Figure \ref{fig3}.
%

We build the feedback control from the value function that was obtained from our procedure and run test examples.
Table \ref{tab:burg} shows the improvement of our method by evaluating the cost functional of the controlled problem over the cost functional of the uncontrolled problem.
We ran the simulation with 10 different initial conditions for each parameter configuration and computed the mean (Column 1), the maximum (Column 2) and the minimum (Column 3) ratio $J^\text{C}/J^\text{UC}$ of the cost functional value for the uncontrolled dynamics $J^\text{UC}$ and controlled dynamics $J^\text{C}$.
We see an improvement over the uncontrolled dynamics for increasing the weights in the cost functional.
It is noteworthy that for almost all parameter configurations and initial values, the LQR controlled dynamics lead to larger cost functional values than the uncontrolled dynamics, which highlights the suboptimality of the LQR control for this scenario.
\begin{table}[htbp]
 \centering
  \begin{tabular}{rrrr}
    \toprule
    & Mean & Best & Worst \\
    \midrule
    $a=0.1$ & 0.985 & 0.969 & 0.996\\
    $a=2.5$ & 0.720 & 0.512 & 0.903\\
     $a=5$ & 0.695 & 0.514 & 0.895\\
    \bottomrule
  \end{tabular}
  \caption{Test 3: The mean, minimum and maximum ratio $J^\text{C}/J^\text{UC}$ of the uncontrolled to the controlled cost for $10$ randomly chosen initial values. Row-wise parameters: $\mu=(a,a)$ with $a = 0.01, 2.5, 5$}
  \label{tab:burg}
\end{table}

\subsection{About the Calculation Times}
Finally, we want to make general remarks about the CPU time, and in particular we discuss the benefit of the offline/online decomposition.
We show in Table \ref{tab:cpu} those results.
The second, third and fourth column refer to the offline costs to compute the basis functions, the VI in the barycenter of each subregion of the parameter space and to precompute the quantities of the affine decomposition, respectively.
We then show the time needed to compute the value function with the PI algorithm. 
Furthermore, we show the benefit of the pre-calculation that speeds up the convergence of the method.
All numbers present average measurements over $10$ parameters drawn randomly from the parameter sets.
We can observe a speed up of factor 2 in the first case, of factor 13 in the second case and of 7 in the third case.
\begin{table}[htbp]
  \centering
  \begin{tabular}{cccccc}
    \toprule
     & \multicolumn{3}{c}{\bf Offline Stage} & \multicolumn{2}{c}{\bf Online Stage}\\
    & Basis gen. & VI  & Precalc. & PI(no precalc.) & PI (precalc.) \\
    \midrule 

 Test 1 & 480  & 72 & 99 &  21 & 11 \\ 
 Test 2 & 12 & 160 & 7 & 161 & 13 \\
 Test 3 & 52 & 25 & 0.3 & 4.5 & 0.6 \\
 \bottomrule
\end{tabular}
\caption{Table for offline and online calculation times (in seconds) for all examples.}
\label{tab:cpu}
\end{table}

\section{Concluding Remarks} \label{sec:conc}
In this paper we have presented an algorithm for the computation of feedback control for parametrized PDEs. Feedback control via DPP suffers from the curse of dimensionality and therefore we make use of model reduction techniques. In the current work, we show that basis functions from the Riccati equation allow to approximate the low-dimensional value function for any parameter. Furthermore, this model reduction approach is very general and allows the computation of the control for any initial conditions since. Finally, we presented an automatic way to generate the reduced domain for the HJB equation. We also want to emphasize that our numerical tests approximate optimal feedback laws for 2D PDEs.

\section*{Acknowledgements}
The authors would like to thank Max Gunzburger for several fruitful discussions
on the subject and valuable comments for improvement of the presentation.

\bibliographystyle{siamplain}
\bibliography{AGHS}

\end{document}